\documentclass[12pt]{article}

\usepackage{amstext}
\usepackage{latexsym}
\usepackage{amssymb}

\def\R{\mathbb R}

\newcommand{\eps}{\varepsilon}

\begin{document}

\newtheorem{theorem}{Theorem}[section]
\renewcommand{\thetheorem}{\arabic{section}.\arabic{theorem}}
\newtheorem{definition}[theorem]{Definition}
\newtheorem{deflem}[theorem]{Definition and Lemma}
\newtheorem{lemma}[theorem]{Lemma}
\newtheorem{example}[theorem]{Example}
\newtheorem{remark}[theorem]{Remark}
\newtheorem{remarks}[theorem]{Remarks}
\newtheorem{cor}[theorem]{Corollary}
\newtheorem{pro}[theorem]{Proposition}
\newtheorem{proposition}[theorem]{Proposition}

\renewcommand{\theequation}{\thesection.\arabic{equation}}

\title{A second look at the Kurth solution in galactic dynamics}
\author{{\sc Markus Kunze$^{1}$} \\[2ex]
      $^{1}$ Universit\"at K\"oln, Institut f\"ur Mathematik, Weyertal 86-90, \\
      D\,-\,50931 K\"oln, Germany \\
      email: mkunze@mi.uni-koeln.de \\[1ex]
      }

\maketitle

\begin{abstract}
\noindent 
The Kurth solution is a particular non-isotropic steady state solution 
to the gravitational Vlasov-Poisson system. It has the property that 
by means of a suitable time-dependent transformation it can be turned 
into a family of time-dependent solutions. Therefore, for a general steady state 
$Q(x, v)=\tilde{Q}(e_Q, \beta)$, depending upon the particle energy $e_Q$ 
and $\beta=\ell^2=|x\wedge v|^2$, the question arises 
if solutions $f$ could be generated that are of the form 
\[ f(t)=\tilde{Q}\Big(e_Q(R(t), P(t), B(t)), B(t)\Big) \] 
for suitable functions $R$, $P$ and $B$, all depending on $(t, r, p_r, \beta)$ 
for $r=|x|$ and $p_r=\frac{x\cdot v}{|x|}$. 
We are going to show that, under some mild assumptions, 
basically if $R$ and $P$ are independent of $\beta$, and if $B=\beta$ is constant, 
then $Q$ already has to be the Kurth solution. 

This paper is dedicated to the memory of Professor Robert Glassey. 
\end{abstract}


\setcounter{equation}{0}

\section{Introduction} 

It is a remarkable fact that very few of Bob Glassey's influential papers 
concern the Vlasov-Poisson system. Certainly there is an in-depth treatment 
of the existence of global solutions in his book \cite{gl1}, 
but apart from that only \cite{gl2,gl3} comes to this author's mind. 
Maybe this is due to Glassey's mathematical formation in the tradition of 
John, Nirenberg, Segal, Strauss ... that he liked better hyperbolic equations, 
and in particular the relativistic Vlasov-Maxwell system. 

For this reason, Glassey would have probably not paid much attention to the present paper, 
but being a polite person, he would nevertheless have found some friendly words for it. 
In addition, this paper has no hard analytic proofs, 
which Glassey could do so well. Let us only mention \cite{gl4} on global existence 
for the `$2.5$' dimensional relativistic Vlasov-Maxwell system, jointly with Jack Schaeffer, 
which is not so well-known (in the sense that not many people have read 
it in all detail), but which is a true masterpiece. One has to use all kinds of structures 
in the system and is not allowed to loose the tiniest part of an $\eps$ 
to close the argument in the end. 
\medskip

Here we are going to consider the Vlasov-Poisson system in the gravitational case, 
which is given by
\begin{equation}\label{vpgr1} 
   \partial_t f(t, x, v)+v\cdot\nabla_x f(t, x, v)-\nabla_x U_f(t, x)\cdot\nabla_v f(t, x, v)=0,
\end{equation} 
where 
\begin{equation}\label{vpgr2} 
   \Delta_x U_f(t, x)=4\pi\rho_f(t, x),\quad\lim_{|x|\to\infty} U_f(t, x)=0,
   \quad\rho_f(t, x)=\int_{\R^3} f(t, x, v)\,dv,
\end{equation} 
for $(t, x, v)\in\R\times\R^3\times\R^3$. Therefore 
\begin{equation}\label{vpgr2b} 
   U_f(t, x)=-\int_{\R^3}\frac{\rho_f(t, y)}{|y-x|}\,dy.
\end{equation} 
The system possesses an abundance of solutions $Q=Q(x, v)$ 
that are independent of time. Let 
\[ e_Q(x, v)=\frac{1}{2}\,|v|^2+U_Q(x) \] 
denote the particle energy and let 
\[ \ell^2(x, v)=|L|^2=|x|^2 |v|^2-(x\cdot v)^2 \] 
be the square of the angular momentum $L=x\wedge v$, respectively. 
Then both $e_Q$ and $\ell^2$ are conserved along solutions 
of the characteristic equations $\ddot{X}(s)=-\nabla U_Q(X(s))$; 
note that also $U_Q$ is independent of time. 
Next recall that a function $g=g(x, v)$ is said to be spherically symmetric, 
if $g(Ax, Av)=g(x, v)$ for all $A\in {\rm SO}(3)$ and $x, v\in\R^3$.
Now it is the content of Jeans's theorem that the distribution function $Q$ 
of every spherically symmetric steady state solution has to be of the form 
\[ Q(x, v)=\tilde{Q}(e_Q(x, v), \ell^2(x, v)) \] 
for a suitable function $\tilde{Q}$ of two variables; see \cite[Section 2]{BFH} 
for a precise formulation. Such steady state solutions are called non-isotropic, 
in contrast to the isotropic ones, which can be written as $Q(x, v)=\tilde{Q}(e_Q(x, v))$; 
a solution of the latter form will necessarily be spherically symmetric, \cite{GNN,R1}. 
\medskip 

In this paper we will have a closer look at 
one particular and non-isotropic steady state solution $Q$, 
which has been found by Kurth in 1978 and which will be denoted 
by $Q_{{\rm K}}$ in the sequel; see \cite{RK}. It is surrounded 
by time-periodic solutions $f_\eps(t)$ such that $f_\eps\to Q_{{\rm K}}$ 
as $\eps\to 0$. Since the $f_\eps(t)$ are semi-explicit, the Kurth solution 
is a good testing ground for all kinds of questions, including some from numerics \cite{RaRe}. 
It is very degenerate in many respects, so an important issue is 
to understand whether it reflects what happens `generically' close to steady states 
(in a sense to be made precise), or on the contrary it is just a peculiarity. 
The Kurth solution is given by 
\begin{eqnarray}\label{Qdef} 
   & & Q_{{\rm K}}(x, v)=\frac{3}{4\pi^3}\,\frac{1}{(1-|x|^2-|v|^2+|x\wedge v|^2)^{1/2}}
   \quad\mbox{where}\quad (\ldots)>0
   \,\,\mbox{and}\,\,|x\wedge v|<1,
   \nonumber
   \\[1ex] & & \mbox{and}\quad Q_{{\rm K}}(x, v)=0\quad\mbox{else}, 
\end{eqnarray} 
for $x, v\in\R^3$. Then (see Lemma \ref{rhokurth} below) its charge density
\[ \rho_{Q_{{\rm K}}}(x)=\int_{\R^3} Q_{{\rm K}}(x, v)\,dv=\frac{3}{4\pi}\,{\bf 1}_{B_1(0)}(x) \] 
is, up to a factor, the characteristic function of the unit ball in $\R^3$. 
The solution to $\Delta U_{Q_{{\rm K}}}=4\pi\rho_{Q_{{\rm K}}}$ 
and $U_{Q_{{\rm K}}}(x)\to 0$ as $|x|\to\infty$ is given by 
\begin{equation}\label{UQkurth} 
   U_{Q_{{\rm K}}}(x)=\left\{\begin{array}{c@{\quad:\quad}c} 
   \frac{1}{2}\,|x|^2-\frac{3}{2} & |x|\le 1
   \\[1ex] -\frac{1}{|x|} & |x|>1
   \end{array}\right. .
\end{equation} 
Next consider the second-order ODE $\ddot{\phi}=-\frac{1}{\phi^2}+\frac{1}{\phi^3}$ 
and let $\phi_\eps$ denote the solution such that $\phi_\eps(0)=1$ and $\dot{\phi}_\eps(0)=\eps$. 
It follows that $\phi_\eps$ is periodic for $\eps<1$ (in fact $|\eps|<1$), and its period is calculated to be
$T_\eps=\frac{2\pi}{(1-\eps^2)^{3/2}}$. Defining
\begin{equation}\label{kurthf} 
   f_\eps(t, x, v)=Q_{{\rm K}}\Big(\frac{x}{\phi_\eps(t)}, \phi_\eps(t)v-\dot{\phi}_\eps(t)x\Big),
   \quad t\in\R,\quad x, v\in\R^3,  
\end{equation}  
the $f_\eps(t)$ are $T_\eps$-periodic and (formal) solutions 
to the gravitational Vlasov-Poisson system; see Lemma \ref{feps_sol} below. 
We may also write 
\begin{equation}\label{faufMQ} 
   f_\eps(t)=Q_{{\rm K}}\circ\Lambda_\eps(t),
   \quad\Lambda_\eps(t)(x, v)=\Big(\frac{x}{\phi_\eps(t)}, \phi_\eps(t)v-\dot{\phi}_\eps(t)x\Big).  
\end{equation}  
The associated density is 
\[ \rho_\eps(t, x)=\int_{\R^3} f_\eps(t, x, v)\,dv=\frac{3}{4\pi}
   \,\frac{1}{\phi_\eps(t)^3}\,{\bf 1}_{\{|x|<\phi_\eps(t)\}}
   =\frac{1}{\phi_\eps(t)^3}\,\rho_{Q_{{\rm K}}}\Big(\frac{x}{\phi_\eps(t)}\Big), \] 
resulting in the potential 
\begin{equation}\label{Ueps} 
   U_\eps(t, x)=\frac{1}{\phi_\eps(t)}\,U_{Q_{{\rm K}}}\Big(\frac{x}{\phi_\eps(t)}\Big).
\end{equation} 
\smallskip

The function $Q_{{\rm K}}$ is spherically symmetric, hence so is $f_\eps(t)$, since  
\[ f_\eps(t, Ax, Av)=Q_{{\rm K}}\Big(A\frac{x}{\phi_\eps(t)}, A[\phi_\eps(t)v-\dot{\phi}_\eps(t)x]\Big)
   =Q_{{\rm K}}\Big(\frac{x}{\phi_\eps(t)}, \phi_\eps(t)v-\dot{\phi}_\eps(t)x\Big)=f_\eps(t, x, v) \] 
for $t\in\R$, $A\in {\rm SO}(3)$, $x, v\in\R^3$. Therefore we may re-express everything 
in the adapted spherically symmetric variables 
\[ r=|x|,\quad p_r=\frac{x\cdot v}{|x|},\quad\ell=|x\wedge v|. \] 
To begin with, 
\begin{eqnarray} 
   e_{Q_K}(r, p_r, \ell^2) & = & \frac{1}{2}\,|v|^2+U_{Q_{{\rm K}}}(r)
   =\frac{1}{2}\,p_r^2+U_{{\rm eff, K}}(r, \ell^2),
   \label{eQkurth} 
   \\[1ex] \quad U_{{\rm eff, K}}(r, \ell^2) & = & U_{Q_{{\rm K}}}(r)+\frac{\ell^2}{2r^2}
   =\left\{\begin{array}{c@{\quad:\quad}c}
   \frac{r^2}{2}-\frac{3}{2}+\frac{\ell^2}{2r^2} & r\le 1 
   \\[1ex] -\frac{1}{r}+\frac{\ell^2}{2r^2} & r\ge 1 \end{array}\right. , 
   \label{classham}
\end{eqnarray}
by (\ref{UQkurth}). Here $U_{{\rm eff, K}}$ is called the effective potential, 
and henceforth we will sometimes write $e$ instead of $e_{Q_{{\rm K}}}$. Also, 
\begin{equation}\label{Nintro} 
   1-|x|^2-|v|^2+|x\wedge v|^2=1-r^2-p_r^2-\frac{\ell^2}{r^2}+\ell^2
   =-2(1+e)+\ell^2,
\end{equation}  
so that 
\begin{eqnarray}\label{QeL} 
   & & \tilde{Q}_{{\rm K}}(e, \ell^2)=\frac{3}{4\pi^3}\,\frac{1}{(-2(1+e)+\ell^2)^{1/2}}
   \quad\mbox{where}\quad (\ldots)>0
   \,\,\mbox{and}\,\,\ell<1,
   \nonumber
   \\[1ex] & & \mbox{and}\quad \tilde{Q}_{{\rm K}}(e, \ell^2)=0\quad\mbox{else}. 
\end{eqnarray}  
In spherically symmetric coordinates, the $\Lambda_\eps(t)$ from (\ref{faufMQ}) are identified with
\begin{equation}\label{got9} 
   \Lambda_\eps(t)(r, p_r)
   =\Big(\frac{r}{\phi_\eps(t)}, \phi_\eps(t)p_r-\dot{\phi}_\eps(t)r\Big),
\end{equation} 
since $f_\eps(t)=Q_{{\rm K}}\circ\Lambda_\eps(t)$ and 
\begin{eqnarray*} 
   \lefteqn{1-\Big|\frac{x}{\phi_\eps(t)}\Big|^2-|\phi_\eps(t)v-\dot{\phi}_\eps(t)x|^2
   +\Big|\frac{x}{\phi_\eps(t)}\wedge [\phi_\eps(t)v-\dot{\phi}_\eps(t)x]\Big|^2}
   \\ & = & 1-\frac{r^2}{\phi_\eps(t)^2}-\phi_\eps(t)^2\,\Big(p_r^2+\frac{\ell^2}{r^2}\Big)
   +2\phi_\eps(t)\dot{\phi}_\eps(t)\,r p_r-\dot{\phi}_\eps(t)^2 r^2+\ell^2
   \\ & = & 1-\frac{r^2}{\phi_\eps(t)^2}-(\phi_\eps(t) p_r-\dot{\phi}_\eps(t) r)^2
   -\phi_\eps(t)^2\frac{\ell^2}{r^2}+\ell^2
   =F\Big(\frac{r}{\phi_\eps(t)}, \phi_\eps(t)p_r-\dot{\phi}_\eps(t)r, \ell^2\Big)
\end{eqnarray*} 
for 
\begin{equation}\label{kiv} 
   F(r, p_r, \ell^2)=1-r^2-p_r^2-\frac{\ell^2}{r^2}+\ell^2.
\end{equation} 
\medskip 

Now we are in position to describe the main result of this paper. 
Writing $\beta=\ell^2$, the Kurth solution is 
\[ Q_{{\rm K}}=\tilde{Q}_{{\rm K}}(e_Q, \beta)
   =\tilde{Q}_{{\rm K}}\Big(e_Q(r, p_r, \beta), \beta\Big), \] 
whereas the neighboring $f_\eps(t)=f_\eps(t, r, p_r, \beta)$ can be expressed as 
\[ f_\eps(t)=\tilde{Q}_{{\rm K}}\Big(e_Q(R_\eps(t), P_\eps(t), B_\eps(t)), B_\eps(t)\Big) \] 
for 
\begin{eqnarray} 
   R_\eps(t) & = & R_\eps(t, r, p_r, \beta)=\frac{r}{\phi_\eps(t)},
   \label{R_kurth} 
   \\[1ex] P_\eps(t) & = & P_\eps(t, r, p_r, \beta)=\phi_\eps(t)p_r-\dot{\phi}_\eps(t)r,
   \label{P_kurth} 
   \\[1ex] B_\eps(t) & = & B_\eps(t, r, p_r, \beta)=\beta,
   \label{B_kurth} 
\end{eqnarray} 
according to (\ref{got9}), (\ref{Nintro}) and (\ref{kiv}). 
It should be remarked (as is verified in Lemma \ref{Zham} below) that defining 
\[ H_{\eps, {\rm K}}(t, r, p_r, \beta)=-\frac{\dot{\phi}_\eps(t)}{\phi_\eps(t)}\,r\,p_r
   -\frac{1}{2}\,(\dot{\phi}_\eps(t)^2-\phi_\eps(t)\ddot{\phi}_\eps(t))\,r^2, \] 
then 
\begin{equation}\label{ham_kurth} 
   \frac{d}{dt}\,Z_\eps(t)=J\nabla H_{\eps, {\rm K}}(t, Z_\eps(t), B_\eps(t)),
\end{equation} 
where 
\[ Z_\eps=(R_\eps, P_\eps),\quad J=\left(\begin{array}{cc} 0 & 1 \\ -1 & 0
   \end{array}\right). \] 
In other words, the time evolution of $Z_\eps$ is governed by 
the time-dependent Hamiltonian $H_{\eps, {\rm K}}$. Also note that both 
$R_\eps$ and $P_\eps$ are in fact independent of $\beta$, and $B_\eps=\beta$ is constant. 
\medskip

Thus, for a general steady state 
$Q(x, v)=\tilde{Q}(e_Q, \beta)$, the question arises 
if solutions $f$ could be found that are of the form 
\begin{equation}\label{f_ast} 
   f(t)=\tilde{Q}\Big(e_Q(R(t), P(t), B(t)), B(t)\Big)
\end{equation}  
for suitable functions $R$, $P$ and $B$, all depending on $(t, r, p_r, \beta)$, 
such that the evolution of $Z=(R, P)$ is Hamiltonian. 
For the moment it will play no role if the $f$ come in a family of $f_\eps$ 
that is close to $Q$ as $\eps\to 0$, or if the function(s) are periodic or not. 
\smallskip

We are going to show that, basically, if $R$ and $P$ are independent of $\beta$, 
and if $B=\beta$ is constant, then $Q$ already has to be the Kurth solution $Q_{{\rm K}}$. 
 
\begin{theorem}\label{mainthm} 
Suppose that the functions 
\[ (R, P, B)(t)=(R(t, r, p_r, \beta), P(t, r, p_r, \beta), B(t, r, p_r, \beta)) \]   
are such that $f(t)=f(t, r, p_r, \beta)$ 
is a solution to the gravitational Vlasov-Poisson system and moreover $R$ and $P$ 
are independent of $\beta$, and $B=\beta$ is constant: 
\[ (R, P, B)(t)=(R(t, r, p_r), P(t, r, p_r), \beta). \] 
Let there exist a Hamiltonian $H=H(t, r, p_r)$ such that 
$\partial_t Z=J\nabla H(t, Z)$ is satisfied for $Z=(R, P)$. 
In addition, we assume that 
\begin{itemize}
\item[(a)] $\partial_e\tilde{Q}\neq 0$ on the support of $Q$ 
and $U'_Q(0)=0$; 
\item[(b)] $\partial_r R>0$ and there is a function $\sigma(t)$ 
such that 
\begin{equation}\label{4farb} 
   \lim_{\delta\to 0}\frac{R(t, \delta, p_r)}{\delta}=\sigma(t);
\end{equation} 
\item[(c)] for the Jacobian of the map $(r, p_r)\mapsto (R(0, r, p_r), P(0, r, p_r))$ 
we have 
\begin{equation}\label{thhue} 
   \det\Big(\frac{\partial (R(0, r, p_r), P(0, r, p_r))}{\partial (r, p_r)}\Big)=1.
\end{equation}  
\end{itemize} 
Then, defining $\alpha=U''_Q(0)$, we must have 
\[ U'_Q(r)=\alpha r,\quad R(t, r)=\frac{r}{\phi(t)},
   \quad P(t, r, p_r)=\phi(t) p_r-\dot{\phi}(t) r, \] 
where $\phi$ solves  
\[ \ddot{\phi}(t)=\alpha\,\Big(-\frac{1}{\phi(t)^2}+\frac{1}{\phi(t)^3}\Big). \] 
\end{theorem} 

\begin{remark}{\rm (a) The proof in Section \ref{prf_sect} is a physics-style 
calculation, and we are not very precise about, for instance, the regularity 
of $\tilde{Q}$. However, this is not the main focus of the paper and missing details 
could be filled in easily. 
\smallskip 

\noindent 
(b) The constant on the right-hand side of (\ref{thhue}) needs not be $1$, 
any other number $\neq 0$ would also work. 
\smallskip 

\noindent 
(c) Concerning hypothesis (b), it will turn out in the proof 
that $\partial_r R\neq 0$, see (\ref{ist1}) below. 
Thus we are going to assume $\partial_r R>0$ without loss of generality. 
Also $R=R(t, r)$ will be shown to be independent of $p_r$. 
Hence (\ref{4farb}) means that in fact 
\[ \sigma(t)=\lim_{\delta\to 0}\frac{R(t, \delta)}{\delta}
   =\partial_r R(t, 0) \] 
is required to exist. This can be guaranteed for instance 
if we suppose that $H\in C^2$. 
\smallskip 

\noindent 
(d) When we started to look into the question if in general solutions 
of the form (\ref{f_ast}) could be found, this author was convinced 
that $\beta$ should play no role for the argument, in the sense that 
everything will be constant in $\beta$. Theorem \ref{mainthm} indicates that 
actually the situation is much more complicated, and that, in a vague sense, 
some `phase mixing' would be needed in order that (\ref{f_ast}) 
could provide a time-dependent solution. 
{\hfill$\diamondsuit$}
}
\end{remark}
\medskip


\setcounter{equation}{0}

\section{Proof of Theorem \ref{mainthm}} 
\label{prf_sect} 

Since the solution $f$ is spherically symmetric by (\ref{f_ast}), 
its potential $U_f$ satisfies
\begin{eqnarray}\label{Ust_lsg}  
   \partial_r U_f (t, \tilde{r}) & = & \frac{4\pi}{\tilde{r}^2}
   \int_0^{\tilde{r}} r^2\rho_f(t, r)\,dr
   =\frac{4\pi}{\tilde{r}^2}\int_0^{\tilde{r}} dr\,r^2\int dv f(t, x, v)
   \nonumber
   \\ & = & \frac{4\pi^2}{\tilde{r}^2}\int_0^{\tilde{r}} dr
   \int\int\,dp_r\,d\beta\,\tilde{Q}\Big(e_Q(R(t), P(t), B(t)), B(t)\Big), 
\end{eqnarray} 
where the arguments of $(R, P, B)$ are $(t, r, p_r, \beta)$ 
and we have used that $dv=\frac{2\pi}{r^2}\,dp_r\,d\ell\,\ell
=\frac{\pi}{r^2}\,dp_r\,d\beta$. 
By definition, 
\[ f(t, r, p_r, \beta)
   =\tilde{Q}\Big(\frac{1}{2}\,P(t)^2+U(R(t))+\frac{B(t)}{2R(t)^2}, B(t)\Big) \] 
for $U=U_Q$ denoting the potential generated by the steady state $Q$. 
From the spherically symmetric version of the Vlasov equation (see \cite{BFH}) 
we hence obtain, for all $(t, r, p_r, \beta)$, 
\begin{eqnarray}\label{ziss} 
   0 & = & \partial_t f(t, r, p_r, \beta)+p_r\,\partial_r f(t, r, p_r, \beta)
   +\Big(\frac{\beta}{r^3}-\partial_r U_f(t, r)\Big)\,\partial_{p_r} f(t, r, p_r, \beta)
   \nonumber
   \\ & = & (\partial_e\tilde{Q})\Big[P(\partial_t P)+U'(R)(\partial_t R)
   +\frac{1}{2R^2}\,(\partial_t B)-\frac{B}{R^3}\,(\partial_t R)\Big]
   +(\partial_\beta\tilde{Q})(\partial_t B)
   \nonumber
   \\ & & +\,p_r (\partial_e\tilde{Q})\Big[P(\partial_r P)+U'(R)(\partial_r R)
   +\frac{1}{2R^2}\,(\partial_r B)-\frac{B}{R^3}\,(\partial_r R)\Big]
   +p_r (\partial_\beta\tilde{Q})(\partial_r B)
   \nonumber
   \\ & & +\,\Big(\frac{\beta}{r^3}-\partial_r U_f(t, r)\Big)
   (\partial_e\tilde{Q})\Big[P(\partial_{p_r} P)+U'(R)(\partial_{p_r} R)
   +\frac{1}{2R^2}\,(\partial_{p_r} B)-\frac{B}{R^3}\,(\partial_{p_r} R)\Big]
   \nonumber
   \\ & & +\,\Big(\frac{\beta}{r^3}-\partial_r U_f(t, r)\Big)
   (\partial_\beta\tilde{Q})(\partial_{p_r} B)
   \nonumber
   \\ & = & (\partial_e\tilde{Q}) P\,
   \Big[\partial_t P+p_r (\partial_r P)
   +\Big(\frac{\beta}{r^3}-\partial_r U_f(t, r)\Big)(\partial_{p_r} P)\Big]
   \nonumber
   \\ & & +(\partial_e\tilde{Q})\Big(U'(R)-\frac{B}{R^3}\Big)
   \,\Big[\partial_t R+p_r (\partial_r R)
   +\Big(\frac{\beta}{r^3}-\partial_r U_f(t, r)\Big)\,(\partial_{p_r} R)\Big]
   \nonumber
   \\ & & +\Big((\partial_e\tilde{Q})\frac{1}{2R^2}+\partial_\beta\tilde{Q}\Big)\,
   \Big[\partial_t B+p_r (\partial_r B)
   +\Big(\frac{\beta}{r^3}-\partial_r U_f(t, r)\Big)
   (\partial_{p_r} B)\Big]. 
\end{eqnarray} 
Since $B=\beta$ is constant by hypothesis, the last line drops out. 
Also $\partial_e\tilde{Q}\neq 0$ on the support of $Q$, 
whence (\ref{ziss}) reduces to 
\begin{eqnarray}\label{fiss}
   0 & = & P\,\Big[\partial_t P+p_r (\partial_r P)
   +\Big(\frac{\beta}{r^3}-\partial_r U_f(t, r)\Big)(\partial_{p_r} P)\Big]
   \nonumber
   \\ & & +\,\Big(U'(R)-\frac{\beta}{R^3}\Big)
   \,\Big[\partial_t R+p_r (\partial_r R)
   +\Big(\frac{\beta}{r^3}-\partial_r U_f(t, r)\Big)\,(\partial_{p_r} R)\Big]. 
\end{eqnarray} 
As $R$ and $P$ are assumed to be independent of $\beta$, 
we can compare the coefficients of the powers $\beta^0, \beta^1, \beta^2$ in $\beta$ 
to deduce that 
\begin{eqnarray} 
   0 & = & P\Big[\partial_t P+p_r (\partial_r P)
   -\partial_r U_f(t, r) (\partial_{p_r} P)\Big] 
   \nonumber
   \\ & & +\,U'(R)\Big[\partial_t R+p_r (\partial_r R)
   -\partial_r U_f(t, r)(\partial_{p_r} R)\Big], 
   \label{strob1} 
   \\ 0 & = & P\frac{1}{r^3}(\partial_{p_r} P)
   +U'(R)\frac{1}{r^3}(\partial_{p_r} R)
   -\frac{1}{R^3}\,\Big[\partial_t R+p_r (\partial_r R)
   -\partial_r U_f(t, r)(\partial_{p_r} R)\Big], 
   \label{strob2} 
   \\ 0 & = & \frac{1}{r^3 R^3}(\partial_{p_r} R). 
   \nonumber
\end{eqnarray} 
Thus $\partial_{p_r} R=0$ and (\ref{strob1}), (\ref{strob2}) 
simplify to  
\begin{eqnarray} 
   0 & = & P\Big[\partial_t P+p_r (\partial_r P)
   -\partial_r U_f(t, r) (\partial_{p_r} P)\Big] 
   +U'(R)\Big[\partial_t R+p_r (\partial_r R)\Big], 
   \label{trae1} 
   \\ 0 & = & P\frac{1}{r^3}(\partial_{p_r} P)
   -\frac{1}{R^3}\,\Big[\partial_t R+p_r (\partial_r R)\Big], 
   \label{trae2} 
\end{eqnarray}
which is a PDE system for $(R(t, r), P(t, r, p_r))$. 
Coming back to (\ref{Ust_lsg}), we have 
\begin{eqnarray}\label{gema}  
   \partial_r U_f (t, \tilde{r}) 
   & = & \frac{4\pi^2}{\tilde{r}^2}\int_0^{\tilde{r}} dr
   \int\int\,dp_r\,d\beta\,
   \tilde{Q}\Big(\frac{1}{2}\,P(t)^2+U(R(t))+\frac{\beta}{2R(t)^2}, \beta\Big)
   \nonumber
   \\ & = & \frac{4\pi^2}{\tilde{r}^2}\int_0^\infty d\beta\int_0^{\tilde{r}} dr\int\,dp_r\,
   \tilde{Q}\Big(\frac{1}{2}\,P(t)^2+U(R(t))+\frac{\beta}{2R(t)^2}, \beta\Big). 
\end{eqnarray}    
Let $\phi_t$ denote the solution map that is associated to 
$\partial_t z=J\nabla H(t, z)$, i.e., $z(t)=\phi_t(r, p_r)$ solves the equation 
and satisfies $z(0)=(r, p_r)$. Thus with $Z(t, r, p_r)=(R(t, r), P(t, r, p_r))$ we get 
\begin{equation}\label{Zvont} 
   Z(t)=\phi_t(Z(0)),
\end{equation}  
since $Z$ is assumed to be a solution. Note that we do not suppose that 
$Z(0)=(r, p_r)$, since this is also not satisfied for the Kurth solution: 
from (\ref{R_kurth}), (\ref{P_kurth}) we have 
$(R_\eps, P_\eps)(0)=(r, p_r-\eps r)$ in this case. 
Since the system is Hamiltonian, each map $\phi_t$ 
is a symplectomorphism \cite[Lemma 1.10]{McDS}, 
and hence in particular $\det D\phi_t=1$ holds for its Jacobian determinant. 
By (\ref{thhue}) from assumption (c) we also know that $\det DZ(0)=1$. 
Therefore (\ref{Zvont}) shows that $\det DZ(t)=1$, and hence 
\[ (\partial_r R)(t, r)(\partial_{p_r} P)(t, r, p_r)
   -(\partial_{p_r} R)(t, r)(\partial_r P)(t, r, p_r)=1, \] 
which in our case is 
\begin{equation}\label{ist1} 
   (\partial_r R)(t, r)(\partial_{p_r} P)(t, r, p_r)=1
\end{equation}
for all $(t, r, p_r)$. At fixed $t$ we are going to apply the change of variables 
\[ (r, p_r)\mapsto (R(t, r), P(t, r, p_r))=Z(t, r, p_r)=(R, P) \] 
to (\ref{gema}), which has $\det DZ(t, r, p_r)=1$ by (\ref{ist1}) 
and $R(t, 0)=0$ due to assumption (b). Then we get 
\begin{eqnarray}\label{glas1} 
   \partial_r U_f (t, \tilde{r}) 
   & = & \frac{4\pi^2}{\tilde{r}^2}\int_0^\infty d\beta\int_0^{R(t, \tilde{r})} dR\int\,dP\,
   \tilde{Q}\Big(\frac{1}{2}\,P^2+U(R)+\frac{\beta}{2R^2}, \beta\Big)
   \nonumber
   \\ & = & \frac{4\pi^2}{\tilde{r}^2}\int_0^{R(t, \tilde{r})} dR\int\,dP\int\,d\beta
   \,\tilde{Q}\Big(\frac{1}{2}\,P^2+U(R)+\frac{\beta}{2R^2}, \beta\Big).  
\end{eqnarray}  
On the other hand, 
\begin{eqnarray}\label{glas2}  
   U'(\tilde{r}) 
   & = & \frac{4\pi}{\tilde{r}^2}\int_0^{\tilde{r}} 
   r^2\rho_Q(r)\,dr
   =\frac{4\pi}{\tilde{r}^2}\int_0^{\tilde{r}} dr\,r^2\int dv\,Q(x, v)
   \nonumber
   \\ & = & \frac{4\pi^2}{\tilde{r}^2}\int_0^{\tilde{r}} dr\int dp_r\int d\beta
   \,\tilde{Q}\Big(\frac{1}{2}\,p_r^2+U(r)+\frac{\beta}{2r^2}, \beta\Big). 
\end{eqnarray} 
Comparing (\ref{glas1}) to (\ref{glas2}), we have shown that  
\begin{equation}\label{umd} 
   \partial_r U_f (t, r)=\frac{R(t, r)^2}{r^2}\,U'(R(t, r))
\end{equation} 
is verified. From (\ref{trae1}), (\ref{umd}) and (\ref{trae2}) 
it follows that
\begin{eqnarray*}
   0 & = & P(\partial_t P+p_r (\partial_r P))
   -P\,\partial_r U_f(t, r) (\partial_{p_r} P)
   +U'(R)(\partial_t R+p_r (\partial_r R))
   \\ & = & P(\partial_t P+p_r (\partial_r P))
   -P\,\frac{R^2}{r^2}\,U'(R)\,(\partial_{p_r} P)
   +U'(R)(\partial_t R+p_r (\partial_r R))
   \\ & = & P(\partial_t P+p_r (\partial_r P))
   -P\,\frac{R^2}{r^2}\,U'(R)\,(\partial_{p_r} P)
   +R^3 P\frac{1}{r^3}(\partial_{p_r} P)\,U'(R)
   \\ & = & \frac{P}{r^3}\Big[r^3(\partial_t P+p_r (\partial_r P))
   -rR^2\,U'(R)\,(\partial_{p_r} P)
   +R^3 (\partial_{p_r} P)\,U'(R)\Big], 
\end{eqnarray*}   
so that   
\begin{equation}\label{carca} 
   r^3(\partial_t P+p_r (\partial_r P))
   +R^2 (R-r)\,U'(R)\,(\partial_{p_r} P)=0.
\end{equation} 
Furthermore, by (\ref{trae2}), 
\[ \partial_{p_r} P^2=2P(\partial_{p_r} P)
   =2\,\frac{r^3}{R^3}\,(\partial_t R+p_r (\partial_r R)), \] 
and since $R$ is independent of $p_r$, integration $\int dp_r$ yields 
\begin{equation}\label{prager} 
   P^2(t, r, p_r)-P^2(t, r,0)
   =2\,\frac{r^3}{R^3}\,p_r\,(\partial_t R)
   +\frac{r^3}{R^3}\,p_r^2\,(\partial_r R).
\end{equation} 
In addition, integration $\int dp_r$ of (\ref{ist1}) leads to 
\[ (\partial_r R)(P(t, r, p_r)-P(t, r, 0))=p_r. \] 
Combining this relation with (\ref{prager}), we get    
\[ P(t, r, p_r)+P(t, r,0)
   =2\,(\partial_r R)\frac{r^3}{R^3}\,(\partial_t R)
   +(\partial_r R)\frac{r^3}{R^3}\,p_r\,(\partial_r R), \] 
which in turn implies that
\begin{eqnarray}\label{witit} 
   P(t, r, p_r) & = & \frac{1}{2}\Big(P(t, r, p_r)-P(t, r, 0)+P(t, r, p_r)+P(t, r, 0)\Big)
   \nonumber
   \\ & = & \frac{p_r}{2\,\partial_r R}
   +(\partial_r R)\,\frac{r^3}{R^3}\,\Big(\partial_t R
   +\frac{1}{2}\,p_r\,(\partial_r R)\Big); 
\end{eqnarray} 
there are only $R$'s on the right-hand side, which are independent of $p_r$.  
If we take the derivative w.r.~to $t$, we obtain 
\begin{eqnarray*}
   \partial_t P & = & -\frac{p_r}{2\,(\partial_r R)^2}\,(\partial^2_{tr} R)
   +(\partial^2_{tr} R)\,\frac{r^3}{R^3}\,\Big(\partial_t R
   +\frac{1}{2}\,p_r\,(\partial_r R)\Big)
   \\ & & -\,(\partial_r R)\,(\partial_t R)\,\frac{3r^3}{R^4}\,
   \Big(\partial_t R+\frac{1}{2}\,p_r\,(\partial_r R)\Big)
   +(\partial_r R)\,\frac{r^3}{R^3}\,\Big(\partial^2_{tt} R
   +\frac{1}{2}\,p_r\,(\partial^2_{tr} R)\Big). 
\end{eqnarray*}  
Due to (\ref{carca}) and (\ref{ist1}), we get  
\begin{eqnarray*} 
   0 & = & r^3\partial_t P+r^3 p_r (\partial_r P)
   +R^2 (R-r)\,U'(R)\,(\partial_{p_r} P)
   \\ & = & r^3\Big[-\frac{p_r}{2\,(\partial_r R)^2}\,(\partial^2_{tr} R)
   +(\partial^2_{tr} R)\,\frac{r^3}{R^3}\,\Big(\partial_t R
   +\frac{1}{2}\,p_r\,(\partial_r R)\Big)
   \\ & & \hspace{1.5em} -\,(\partial_r R)\,(\partial_t R)\,\frac{3r^3}{R^4}\,
   \Big(\partial_t R+\frac{1}{2}\,p_r\,(\partial_r R)\Big)
   +(\partial_r R)\,\frac{r^3}{R^3}\,\Big(\partial^2_{tt} R
   +\frac{1}{2}\,p_r\,(\partial^2_{tr} R)\Big)+p_r (\partial_r P)\Big]
   \\ & & +R^2 (R-r)\,U'(R)\,\frac{1}{(\partial_r R)}. 
\end{eqnarray*}  
Taking $p_r=0$, this yields   
\begin{equation}\label{haberl} 
   0=\frac{r^6}{R^3}
   \Big[(\partial^2_{tr} R)\,(\partial_t R)
   -(\partial_r R)\,(\partial_t R)^2\,\frac{3}{R}
   +(\partial_r R)\,(\partial^2_{tt} R)\Big]
   +R^2 (R-r)\,U'(R)\,\frac{1}{(\partial_r R)}.
\end{equation}   
If $t$ is fixed, then the map $(r, p_r)\mapsto (R(t, r), P(t, r, p_r))=Z(t, r, p_r)$ 
is symplectic, owing to (\ref{ist1}), and it is generated by the `point transformation' 
$r\mapsto R(t, r)$. Thus, using \cite[equ.~(1.44)]{MZ}, there is a scalar function 
$v=v(t, r)$ such that 
\begin{equation}\label{Pmitv} 
   P(t, r, p_r)=\frac{1}{\partial_r R(t, r)}\,(p_r-\partial_r v(t, r)).
\end{equation}  
Therefore, from (\ref{witit}),   
\begin{equation}\label{rk47} 
   \frac{1}{\partial_r R}\,(p_r-\partial_r v(t, r))
   =P=\frac{p_r}{2\,\partial_r R}
   +(\partial_r R)\,\frac{r^3}{R^3}\,\Big(\partial_t R
   +\frac{1}{2}\,p_r\,(\partial_r R)\Big).
\end{equation} 
Taking once again $p_r=0$, we see that   
\begin{equation}\label{drv} 
   \partial_r v(t, r)
   =-(\partial_r R)^2\,(\partial_t R)\,\frac{r^3}{R^3}.
\end{equation}  
If we plug this relation back to (\ref{rk47}), it follows that 
\[ (\partial_r R)^3\,\frac{r^3}{R^3}=1, \] 
or 
\[ \partial_r R=\frac{R}{r}. \] 
For $\delta>0$ integration yields 
\[ R(t, r)=r\,\frac{R(t, \delta)}{\delta}\,e^{C(t)} \] 
for a suitable function $C(t)$. By assumption, taking the limit $\delta\to 0$, 
we get 
\[ R(t, r)=a(t)r, \] 
where $a(t)=\sigma(t)e^{C(t)}$. Thus (\ref{haberl}) simplifies to 
\[ 0=\frac{r^3}{a(t)^3}
   \Big[\dot{a}(t)^2 r-a(t)\,\dot{a}(t)^2 r^2\,\frac{3}{a(t)r}
   +a(t)\ddot{a}(t)r\Big]
   +a(t)^2 r^3 (a(t)-1)\,U'(a(t)r)\,\frac{1}{a(t)}, \] 
which is 
\[ U'(a(t)r)
   =\frac{2\dot{a}(t)^2-a(t)\ddot{a}(t)}{a(t)^4 (a(t)-1)}\,r. \]
Replacing $a(t)r$ by $r$, this leads to 
\[ U'(r)
   =\frac{2\dot{a}(t)^2-a(t)\ddot{a}(t)}{a(t)^5 (a(t)-1)}\,r. \]
Since the variables are separated, we deduce that there is $\alpha\in\R$ such that 
\begin{equation}\label{sepi} 
   \frac{U'(r)}{r}=\alpha
   =\frac{2\dot{a}(t)^2-a(t)\ddot{a}(t)}{a(t)^5 (a(t)-1)}
\end{equation} 
for all $(t, r)$. Thus if we set $\phi(t)=\frac{1}{a(t)}$, 
then $\dot{\phi}(t)=-\frac{\dot{a}(t)}{a(t)^2}$ and, by (\ref{sepi}),  
\[ \ddot{\phi}(t)=-\frac{\ddot{a}(t)}{a(t)^2}+2\,\frac{\dot{a}(t)^2}{a(t)^3}
   =\frac{1}{a(t)^3}\,(2\dot{a}(t)^2-a(t)\ddot{a}(t))
   =\alpha\,a(t)^2 (a(t)-1)=\alpha\,\Big(-\frac{1}{\phi(t)^2}+\frac{1}{\phi(t)^3}\Big). \] 
Moreover, we also have
\[ R(t, r)=\frac{r}{\phi(t)}. \]  
Next, using (\ref{drv}), 
\[ \partial_r v(t, r)
   =-a(t)^2\,\dot{a}(t)r\,\frac{r^3}{a(t)^3 r^3}
   =-\frac{\dot{a}(t)}{a(t)}\,r=\frac{\dot{\phi}(t)}{\phi(t)}\,r, \] 
and hence upon integration 
\[ v(t, r)=\frac{1}{2}\,\frac{\dot{\phi}(t)}{\phi(t)}\,r^2+\gamma(t) \] 
for a suitable function $\gamma(t)$. Thus (\ref{Pmitv}) implies that 
\[ P(t, r, p_r)=\frac{1}{\partial_r R(t, r)}\,(p_r-\partial_r v(t, r)) 
   =\phi(t)\,\Big(p_r-\frac{\dot{\phi}(t)}{\phi(t)}\,r\Big)
   =\phi(t)\,p_r-\dot{\phi}(t)\,r. \] 
Lastly, $U'(0)=0$ in conjunction with (\ref{sepi}) shows that 
\[ \alpha=\lim_{r\to 0}\frac{U'(r)}{r}=U''(0).  \] 
This completes the proof of Theorem \ref{mainthm}. 
{\hfill$\Box$}\bigskip 


\setcounter{equation}{0}

\section{Some technical results} 

\begin{lemma}\label{rhokurth} One has 
\[ \int_{\R^3} Q_{{\rm K}}(x, v)\,dv=\frac{3}{4\pi}\,{\bf 1}_{B_1(0)}(x) \] 
\end{lemma} 
{\bf Proof\,:} We only consider $r=|x|<1$. First note that 
\[ 1-|x|^2-|v|^2+|x\wedge v|^2=1-r^2-p_r^2-\frac{\ell^2}{r^2}+\ell^2>0 \] 
means that $p_r^2<(1-r^2)(1-\frac{\ell^2}{r^2})\le 1-r^2$ and also 
\[ \ell^2<r^2-\frac{r^2 p_r^2}{1-r^2}=:\ell_0^2. \]
In spherical symmetry we have $dv=\frac{2\pi}{r^2}\,dp_r\,d\ell\,\ell$. 
Therefore by (\ref{Qdef}), 
\begin{eqnarray*} 
   \int_{\R^3} Q_{{\rm K}}(x, v)\,dv
   & = & \frac{2\pi}{r^2}\int_{|p_r|\le\sqrt{1-r^2}} dp_r
   \int_0^{\ell_0} d\ell\,\ell\,Q_{{\rm K}}(r, p_r, \ell)
   \\ & = & \frac{3}{2\pi^2 r^2}\int_{|p_r|\le\sqrt{1-r^2}} dp_r
   \int_0^{\ell_0} d\ell\,\ell\,
   \frac{1}{(1-r^2-p_r^2-\frac{\ell^2}{r^2}+\ell^2)^{1/2}}
   \\ & = & \frac{3}{2\pi^2}\,\frac{1}{1-r^2}
   \int_{|p_r|\le\sqrt{1-r^2}} dp_r\int_0^{\ell_0} d\ell\,
   (-1)\frac{d}{d\ell}\,\Big(1-r^2-p_r^2-\frac{\ell^2}{r^2}+\ell^2\Big)^{1/2}
   \\ & = & \frac{3}{2\pi^2}\,\frac{1}{1-r^2}
   \int_{|p_r|\le\sqrt{1-r^2}} dp_r\,(1-r^2-p_r^2)^{1/2}
   \\ & = & \frac{3}{\pi^2}\,\int_0^1 ds\,(1-s^2)^{1/2}
   \\ & = & \frac{3}{4\pi}, 
\end{eqnarray*} 
as was to be shown. 
{\hfill$\Box$}\bigskip 

\begin{lemma}\label{feps_sol} 
The function $f_\eps$ from (\ref{kurthf}) is a (formal) solution 
to the gravitational Vlasov-Poisson system. 
\end{lemma} 
{\bf Proof\,:} To begin with, (\ref{kurthf}) yields 
\begin{equation}\label{dtfeps} 
   \partial_t f_\eps
   =-\frac{\dot{\phi}_\eps}{\phi_\eps^2}\,x\cdot ((\nabla_x Q_{{\rm K}})\circ \Lambda_\eps)
   -\ddot{\phi}_\eps\,x\cdot ((\nabla_v Q_{{\rm K}})\circ \Lambda_\eps)
   +\dot{\phi}_\eps\,v\cdot ((\nabla_v Q_{{\rm K}})\circ \Lambda_\eps),
\end{equation} 
where $\Lambda_\eps$ is defined in (\ref{faufMQ}). 
Observe that on $\{Q_{{\rm K}}\neq 0\}$: 
\begin{eqnarray}
   \nabla_x Q_{{\rm K}}(x, v) & = & \Big(\frac{4\pi^3}{3}\Big)^2 Q_{{\rm K}}^3(x, v)\,((1-|v|^2)x
   +\langle x, v\rangle v),
   \label{nabxkur} 
   \\ \nabla_v Q_{{\rm K}}(x, v) & = & \Big(\frac{4\pi^3}{3}\Big)^2 Q_{{\rm K}}^3(x, v)\,((1-|x|^2)v
   +\langle x, v\rangle x).
   \label{nabvkur}  
\end{eqnarray} 
Furthermore, 
\begin{eqnarray*} 
   \nabla_x f_\eps & = & \frac{1}{\phi_\eps}\,((\nabla_x Q_{{\rm K}})\circ \Lambda_\eps)
   -\dot{\phi}_\eps\,((\nabla_v Q_{{\rm K}})\circ \Lambda_\eps),
   \\ \nabla_v f_\eps & = & \phi_\eps\,((\nabla_v Q_{{\rm K}})\circ \Lambda_\eps). 
\end{eqnarray*} 
On $\{f_\eps\neq 0\}$ we have 
\[ U_\eps(t, x)=\frac{1}{\phi_\eps(t)}\,\bigg(\frac{|x|^2}{2\phi_\eps(t)^2}-\frac{3}{2}\bigg),
   \quad\nabla_x U_\eps(t, x)=\frac{x}{\phi_\eps(t)^3}. \] 
Thus we obtain from $\ddot{\phi}_\eps=-\frac{1}{\phi_\eps^2}+\frac{1}{\phi_\eps^3}$ that 
\begin{eqnarray*} 
   \lefteqn{\partial_t f_\eps+v\cdot\nabla_x f_\eps-\nabla_x U_\eps\cdot\nabla_v f_\eps}
   \\ & = & -\frac{\dot{\phi}_\eps}{\phi_\eps^2}\,x\cdot ((\nabla_x Q_{{\rm K}})\circ \Lambda_\eps)
   -\ddot{\phi}_\eps\,x\cdot ((\nabla_v Q_{{\rm K}})\circ \Lambda_\eps)
   +\dot{\phi}_\eps\,v\cdot ((\nabla_v Q_{{\rm K}})\circ \Lambda_\eps) 
   \\ & & +\,\frac{1}{\phi_\eps}\,v\cdot ((\nabla_x Q_{{\rm K}})\circ \Lambda_\eps)
   -\dot{\phi}_\eps\,v\cdot ((\nabla_v Q_{{\rm K}})\circ \Lambda_\eps)
   -\frac{x}{\phi_\eps^2}\cdot ((\nabla_v Q_{{\rm K}})\circ \Lambda_\eps)
   \\ & = & -\frac{\dot{\phi}_\eps}{\phi_\eps^2}\,x\cdot ((\nabla_x Q_{{\rm K}})\circ \Lambda_\eps)
   -\frac{1}{\phi_\eps^3}\,x\cdot ((\nabla_v Q_{{\rm K}})\circ \Lambda_\eps)
   +\,\frac{1}{\phi_\eps}\,v\cdot ((\nabla_x Q_{{\rm K}})\circ \Lambda_\eps)
   \\ & = & \frac{1}{\phi_\eps^3}\Big(\frac{4\pi^3}{3}\Big)^2 (Q_{{\rm K}}^3\circ\Lambda_\eps)
   \bigg[-\phi_\eps\dot{\phi}_\eps\,x\cdot\Big((1-|\phi_\eps v-\dot{\phi}_\eps x|^2)\frac{x}{\phi_\eps}
   +\langle\frac{x}{\phi_\eps}, \phi_\eps v-\dot{\phi}_\eps x\rangle [\phi_\eps v-\dot{\phi}_\eps x]\Big)
   \\ & & \hspace{10em} -\,x\cdot\Big((1-\Big|\frac{x}{\phi_\eps}\Big|^2)
   [\phi_\eps v-\dot{\phi}_\eps x]+\langle\frac{x}{\phi_\eps}, \phi_\eps v
   -\dot{\phi_\eps}x\rangle\frac{x}{\phi_\eps}\Big)
   \\ & & \hspace{10em} +\,\phi_\eps^2\,v\cdot\Big((1-|\phi_\eps v-\dot{\phi}_\eps x|^2)\frac{x}{\phi_\eps}
   +\langle\frac{x}{\phi_\eps}, \phi_\eps v-\dot{\phi}_\eps x\rangle 
   [\phi_\eps v-\dot{\phi}_\eps x]\Big)\bigg]
   \\ & = & 0, 
\end{eqnarray*} 
where the last step requires some calculation. Apart from that, 
\[ \Delta U_\eps(t, x)=\frac{1}{\phi_\eps(t)^3}\,\Delta U_{Q_{{\rm K}}}\Big(\frac{x}{\phi_\eps(t)}\Big)
   =\frac{4\pi}{\phi_\eps(t)^3}\,\rho_{Q_{{\rm K}}}\Big(\frac{x}{\phi_\eps(t)}\Big)=4\pi\rho_\eps(t, x), \] 
which completes the somewhat formal argument. 
{\hfill$\Box$}\bigskip  

\begin{lemma}\label{Zham} The functions $R_\eps$, $P_\eps$ and $B_\eps$ 
from (\ref{R_kurth}), (\ref{P_kurth}) and (\ref{B_kurth}), respectively, 
provide a solution to (\ref{ham_kurth}). 
\end{lemma} 
{\bf Proof\,:} Since $H_{\eps, {\rm K}}$ is independent of $\beta$, we drop this variable. Then 
\[ \partial_r H_{\eps, {\rm K}}(t, r, p_r)=-\frac{\dot{\phi}_\eps(t)}{\phi_\eps(t)}\,p_r
   -(\dot{\phi}_\eps^2(t)-\phi_\eps(t)\ddot{\phi}_\eps(t))\,r, 
   \quad\partial_{p_r} H_{\eps, {\rm K}}(t, r, p_r)=-\frac{\dot{\phi}_\eps(t)}{\phi_\eps(t)}\,r. \] 
This yields 
\[ \dot{R}_\eps(t)=-\frac{r}{\phi_\eps^2(t)}\,\dot{\phi}_\eps(t)
   =\partial_{p_r} H_{\eps, {\rm K}}(t, R_\eps(t), P_\eps(t)) \] 
as well as 
\begin{eqnarray*} 
   \dot{P}_\eps(t) 
   & = & \dot{\phi}_\eps(t)p_r-\ddot{\phi}_\eps(t)r
   \\ & = & \frac{\dot{\phi}_\eps(t)}{\phi_\eps(t)}\,
   \Big(\phi_\eps(t)p_r-\dot{\phi}_\eps(t)r\Big)
   +(\dot{\phi}_\eps^2(t)-\phi_\eps(t)\ddot{\phi}_\eps(t))\,\frac{r}{\phi_\eps(t)}
   \\ & = & \frac{\dot{\phi}_\eps(t)}{\phi_\eps(t)}\,P_\eps(t)
   +(\dot{\phi}_\eps^2(t)-\phi_\eps(t)\ddot{\phi}_\eps(t))\,R_\eps(t)
   \\ & = & -\partial_r H_{\eps, {\rm K}}(t, R_\eps(t), P_\eps(t)), 
\end{eqnarray*} 
and altogether this yields (\ref{ham_kurth}). 
{\hfill$\Box$}\bigskip  


\end{document}